\documentclass[a4paper]{ifacconf}
\usepackage{cite}
\usepackage{amsmath,amssymb,amsfonts}
\usepackage{algorithmic}
\usepackage{graphicx}
\usepackage{textcomp}
\usepackage{xcolor}
\usepackage{bm}
\usepackage[font=scriptsize]{caption}

\makeatletter
\let\old@ssect\@ssect 
\makeatother

\usepackage{natbib}        
\usepackage{url}
\usepackage[pdfstartview=FitV, urlcolor=blue]{hyperref}

\makeatletter
\def\@ssect#1#2#3#4#5#6{%
  \NR@gettitle{#6}
  \old@ssect{#1}{#2}{#3}{#4}{#5}{#6}
}
\makeatother

\newcommand*{\CopyCounter}[2]{%
  \expandafter\def\csname c@#2\endcsname{\csname c@#1\endcsname}%
  \expandafter\def\csname p@#2\endcsname{\csname p@#1\endcsname}%
  \expandafter\def\csname the#2\endcsname{\csname the#1\endcsname}}

\CopyCounter{Theorem}{ProposedProblem}
\CopyCounter{Theorem}{Proposition}
\CopyCounter{Theorem}{Property}
\CopyCounter{Theorem}{Claim}
\CopyCounter{Theorem}{Lemma}
\CopyCounter{Theorem}{Corollary}
\CopyCounter{Theorem}{Conjecture}
\CopyCounter{Theorem}{Definition}
\CopyCounter{Theorem}{Example}
\CopyCounter{Theorem}{Remark}
\CopyCounter{Theorem}{Question}
\CopyCounter{Theorem}{Condition}
\CopyCounter{Theorem}{Criterion}
\CopyCounter{Theorem}{Observation}
\theoremstyle{plain}

\newcommand{\x}{\bm{x}}

\newcommand{\ba}{\bm{a}}

\newcommand{\Bf}{\bm{f}}

\newcommand{\y}{\bm{y}}

\newcommand{\J}{{\cal{J}}}

\newcommand{\R}{\mathbb{R}}

\newcommand{\vbar}{\bar{v}}

\newcommand{\TG}{T_{\mbox{\tiny{G}}}^{}}
\newcommand{\subG}{_{\mbox{\tiny{G}}}^{}}
\newcommand{\TY}{T_{\mbox{\tiny{Y}}}^{}}
\newcommand{\TR}{T_{\mbox{\tiny{R}}}^{}}
\newcommand{\sg}{s_{\mbox{\tiny{G}}}^{}}
\newcommand{\sr}{s_{\mbox{\tiny{R}}}^{}}
\newcommand{\srsq}{s_{\mbox{\tiny{R}}}^{2}}
\newcommand{\sbeta}{s_{\beta}^{}}
\newcommand{\salpha}{s_{\alpha}^{}}
\newcommand{\shatalpha}{\hat{s}_{\alpha}^{}}
\newcommand{\shatbeta}{\hat{s}_{\beta}^{}}
\newcommand{\subR}{_{\mbox{\tiny{R}}}^{}}
\newcommand{\subY}{_{\mbox{\tiny{Y}}}^{}}

\newcommand{\vc}{v_c^{}}
\newcommand{\dl}{d_{\ell}^{}}
\newcommand{\dalpha}{d_{\alpha}^{}}
\newcommand{\dbeta}{d_{\beta}^{}}
\newcommand{\ug}{q}
\newcommand{\Jg}{\mathcal{J}_{\mbox{\tiny{G}}}^{}}
\newcommand{\wTY}{\widetilde{T}\subY}
\newcommand{\DY}{D\subY}
\newcommand{\DR}{D\subR}

\newcommand{\tinySubs}[1]{_{\mbox{\tiny{#1}}}^{}}

\newcommand{\PP}{\mathbb{P}}

\begin{document}
\begin{frontmatter}

\title{Optimal Driving
Under Traffic Signal Uncertainty\thanksref{footnoteinfo}} 

\thanks[footnoteinfo]{Supported by the NSF DMS (awards 1645643, 1738010, 2111522) and the first author's NDSEG Fellowship. This work has been submitted to IFAC for possible publication.}

\author[First]{Mallory E. Gaspard,  Alexander Vladimirsky} 

\address[First]{Cornell University, 
  Ithaca, NY 14853 USA\\ (e-mail: \{ meg375, vladimirsky \} @cornell.edu).}

\begin{abstract}     
We study driver's optimal trajectory planning under uncertainty in the duration of a traffic light's green phase.
We interpret this as an optimal control problem with an objective of minimizing the expected cost based on the fuel use, discomfort from rapid velocity changes, and time to destination.
Treating this in the framework of dynamic programming, we show that the probability distribution on green phase durations gives rise to a sequence of Hamilton-Jacobi-Bellman PDEs,
which are then solved numerically to obtain optimal acceleration/braking policy in feedback form.
Our numerical examples illustrate the approach and highlight the role of conflicting goals and uncertainty in shaping drivers' behavior.
\end{abstract}

\begin{keyword}
Optimal control under uncertainty; dynamic programming; 
traffic signal advisory.
\end{keyword}
\end{frontmatter}

\section{INTRODUCTION}

Determining optimal driving behaviors under traffic signal uncertainty is of crucial importance in modern transportation applications.  
Recent technological advances aim to provide real-time ``signal phase and timing'' (SPaT) information 
to smartphones and on-board vehicle systems, with the goal of improving safety, fuel efficiency, and driver's comfort.
While this technology will yield likely decreases in crash rates and carbon emissions (\cite{chang2015estimated}), the full transition to vehicles communicating directly with traffic signals will take several decades (\cite{motavalli2020smartlights}), and some level of uncertainty in SPaT information will remain even then due to occasional network failures and the presence of actuated or semi-actuated traffic signals (\cite{eteifa2021predicting}). 

Modeling how drivers do (or should) utilize SPaT knowledge is a popular research area; e.g., \cite{sun2018robust}, \cite{sun2020optimal}, \cite{mahler2014optimal}, \cite{mahler2012reducing}, \cite{zhao2021online}. But typical prior approaches are based on either full SPaT information or ad hoc models of SPaT-robustness, with the emphasis placed on detailed vehicle dynamics, road models with multiple signalized intersections, traffic effects, and route-planning choices to come up with  trip plans that optimize the fuel efficiency. 
In contrast, we consider the simplest scenario (one car with simplified dynamics on an empty road with a single traffic signal -- as described in Section \ref{section:problem_formulation_paper}), focusing instead on a detailed model of SPaT-uncertainty due to the probability distribution over possible durations of the green light (Section \ref{section:uncertainty_models_paper}). 
The driver knows the yellow light and red light durations, so they can deterministically plan the rest of their strategy once the initial green light turns yellow. They must also obey the traffic laws (red light and speed limit) and balance expected costs related to fuel usage, aggressive acceleration / deceleration behavior, and total time to target.

We address this optimal control problem in the dynamic programming (DP) 
framework, as was recently done for a related {\em red phase duration uncertainty} problem in (\cite{typaldos2020vehicle}).
However, the dilemma in the green phase (i.e., should we try to beat the red light?)
is far more affected by the available SPaT information and the duration of the transitional yellow phase.
Unlike in (\cite{typaldos2020vehicle}), where the discrete nature of DP was essential and the timestep could not be refined,
our continuous DP yields a sequence of Hamilton-Jacobi-Bellman (HJB) equations, which are then solved numerically (Section \ref{section:numerical_implementation_paper}) to determine the optimal acceleration and braking policy in feedback form. 
Our numerical experiments in Section \ref{section:numerical_exp_paper} highlight the role of conflicting optimization goals and SPaT-uncertainty patterns in shaping drivers' behavior.

\section{PROBLEM FORMULATION}\label{section:problem_formulation_paper}

Throughout this paper, we focus on the problem of a single driver aiming to optimize their route to a target through an intersection with one traffic light.
We assume the vehicle is just a point mass traveling right-to-left on a flat single-lane road (an interval $[d^*, \bar{d}]$) with no other traffic:
from its starting position at $d \leq \bar{d}$ with initial speed $v$, through the intersection at $\dl < d$ and towards its target at $d^* < \dl.$
The driver varies the car's signed acceleration and must stay under the speed limit $\bar{v}$ on the way to $d^*$.
The car may not enter the intersection (crossing $\dl$) while the light is red, but may do so while the light is yellow\footnote{
This corresponds to {\em permissive} yellow light laws (\cite{beeber2020kinematiceqn}),
where the car does not have to clear the intersection by the time the light turns red.
This makes the actual intersection width and the physical car length irrelevant.}. 
We further assume that the driver knows the yellow and red signal durations $(\DY, \DR).$

{\em Time Horizon. \hspace{1mm}} 
The entire planning process is split into 
(1)  an initial green light phase $I$ starting at $t=0$ followed by 
(2) a yellow phase $Y$  starting at $t=\TY$, 
(3) a red phase $R$ starting at $\TR = \TY + \DY$, and 
(4) an indefinite green phase $G$ starting at $\TG = \TR + \DR$ and terminating when we reach the target at $d^*$.
In many 
realistic situations, the
traffic signal uncertainty will primarily rest in phase $I$. The limited 
knowledge scenario discussed in Section \ref{section:uncertainty_models_paper} will define a probability distribution on the {\em random turning yellow time} $\wTY$, thereby leading to uncertainty in the finite time horizon 
of phase $I$.

{\em Vehicle Dynamics. \hspace{1mm}}
Let $\Omega = \{(d,v,r) : d \in [d^*, \bar{d}], v\in [0, \bar{v}], r \geq 0\}$. 
We will 
use $\x = \begin{bmatrix}
d & v
\end{bmatrix}^T$
to refer to the vehicle's initial configuration or a generic point
in $(d,v)$ space,
switching to  
$\y(r) = \begin{bmatrix}
d(r) & v(r)
\end{bmatrix}^T$ to encode how 
position and velocity change with time.
Starting from an initial configuration $\y(t) = \x$ and assuming
 that external forces acting on the car are negligible, the vehicle's dynamics
are  
\begin{equation}\label{eqn:frictionless_dynamics}
\dot{\y}(r) \, = \, \boldsymbol{f}\left(\y(r), a(r) \right)
\, = \,  
\begin{bmatrix}
-v(r) \\
a(r)
\end{bmatrix}, \qquad \text{for } r \geq t
\end{equation}
where 
the vehicle's acceleration 
$a(\cdot)$ is a measurable control function, $a: \R \to A,$ and the set of available (signed) acceleration values is 
$A = [-\alpha, \, \beta].$

{\em Cost Function. \hspace{1mm}}
The driver's goal on the trip from $d$ to $d^*$ is to select 
$a(\cdot)$ to
address several
optimization objectives: fuel consumption\footnote{
We approximate fuel consumption as the amount of fuel let into the engine. We assume that fuel enters the engine at a rate proportional to the car's positive acceleration (i.e., how far down the driver is pressing the gas pedal) (\cite{bennett2016medium}).}, discomfort from rapid acceleration/deceleration, and time to target. 
We aim to minimize an integral of the {\em running cost}
\begin{equation}\label{eqn:running_cost_ob_paper}
    K(\y(r), a) = c_1 [a]_+ + c_2 a^2 + c_3
\end{equation}
where $[a]_+ = \max\{a, 0\}.$ The non-negative objective weights $c_i$'s
reflect individual preferences.  In principle, these can also be learned from trends in driver behavior data, as in (\cite{butakov2016personalized}).

{\em Deterministic Control Problem. \hspace{1mm}}
We begin by discussing the fully deterministic problem that starts after the light turns yellow  (i.e., the optimal driving for $t \geq \TY$) .
In the usual style of dynamic programming, we consider the remaining phases of the traffic light backwards in time.

We first consider the last green phase of unlimited duration, $G$, which is best described as an {\em exit-time problem},
terminating when the vehicle reaches its target at the time $T^*= \min\{r \geq 
\TG \mid \y(r) = d^*\}.$
The cost-to-go depends on our initial configuration  
and the chosen control $\ba(\cdot)$, 
which together determine the remaining time to target $(T^* - t).$ 
We thus define the control-specific cost 
\begin{equation}\label{eqn:cost_functional}
    \Jg  (\x, \, a(\cdot)) \; = \; \int^{T^*}_t K \left(\y(r), a(r) \right) \, dr,
\end{equation}
starting from $\x = \y(t)$ with $t \geq \TG.$  
By a standard argument (\cite{bardi1997optimalcontrol}), the value function defined as
\begin{equation}\label{eqn:oc_problem_stationary} 
   \ug(\x) =  \inf_{a(\cdot)} \{ \Jg(\x, a(\cdot)) \} 
\end{equation}
can be recovered as a viscosity solution of a 
stationary Hamilton-Jacobi-Bellman (HJB) equation 
\begin{equation}\label{eqn:hjb_stationary}
    0 = \min_{a \in A} \{K(\x, a) + \nabla \ug \cdot \boldsymbol{f}(\x, a)\}
\end{equation}
with boundary conditions 
$
    \ug(d^*, v) = 0,
    \, \forall v \in [0, \vbar].
$

During the other two deterministic stages ($Y$ and $R$), we define the control-specific cost as 
\begin{equation}\label{eqn:cost_functional_yr}
    \mathcal{J}(\x, t, a(\cdot)) = \int^{\TG}_t K \left(\y(r), a(r) \right) dr + \ug(\y(\TG)),
\end{equation} 
starting
from $\x$ at the time $t \in [\TY, \, \TG).$
The value function  
\begin{equation}\label{eqn:valuefn_timedep}
       u(\x, t) =  \inf_{a(\cdot)} \{ \mathcal{J}(\x, t, a(\cdot)) \}
\end{equation}
is now time-dependent to account for the fixed time-horizon.
The state constraints
detailed below
can often 
make this value function discontinuous during the $Y$ phase.
Nevertheless, $u$ can still be recovered as a discontinuous viscosity solution  \cite[Chapter 5]{bardi1997optimalcontrol} of the time-dependent HJB equation
\begin{equation}\label{eqn:hjb_timedep}
    0 = u_t + \min_{a \in A} \{K(\x, a) + \nabla u \cdot \boldsymbol{f}(\x, a) \}
\end{equation}
subject to the terminal condition $u(\x,\TG) = \ug(\x).$  

{\em State Constraints. \hspace{1mm}}
The car should never enter the intersection during the $R$ phase\footnote{
For the sake of brevity,
we omit detailed algebraic derivations in favor of 
geometric discussion emphasizing the physical intuition and main results. 
Further details of the numerical implementation are discussed in Section \ref{section:numerical_implementation_paper}.
}.  This makes it natural to consider a space-time obstacle
$\mathcal{I} = \{(d,v,t) \, \mid \, d = \dl, \, 0 \leq v \leq \bar{v}, \, \TR \leq t < \TG\}, $ but we go further and forbid the car from
taking any $(d,v)$ 
from which it is impossible to come to a complete stop before reaching $\dl.$
(This guarantees that the car will not violate the rules even if the signal fails to switch to green light at the usual time $\TG$.)
The minimal stopping distance depends on the 
current velocity $v$ and the maximum deceleration rate $\alpha$, yielding
the curve of critical positions  (the ``parabola of last resort'') $\dalpha(v) = \dl + \frac{v^2}{2 \alpha}$ and the disallowed region 
$$\mathcal{I}\subR (t) = \{(d,v) \mid  \dl \leq d < \dalpha(v), \, 0 \leq v \leq \bar{v} \}$$
enforced for $t \in [\TR, \, \TG)$.
Since we are forced to use the maximum deceleration $a = - \alpha$ on the parabola of last resort, this makes it easy to compute the 
value function on it.  Let $\sg= \TG - t$. It would take the vehicle $\salpha = v / \alpha$ to come to a full stop, but the actual braking time might be shorter,
$\shatalpha= \min\{s_{\alpha}, \sg\}.$
The value function along $(\dalpha(v), v, t)$ is 
\begin{equation}
\label{eq:c_alpha}
C_{\alpha} = \shatalpha c_2 \alpha^2+ c_3 \sg + q\left(\y(\TG)\right).
\end{equation}

\begin{figure}[ht]
\begin{minipage}[b]{0.45\linewidth}
\centering
\includegraphics[width=\textwidth]{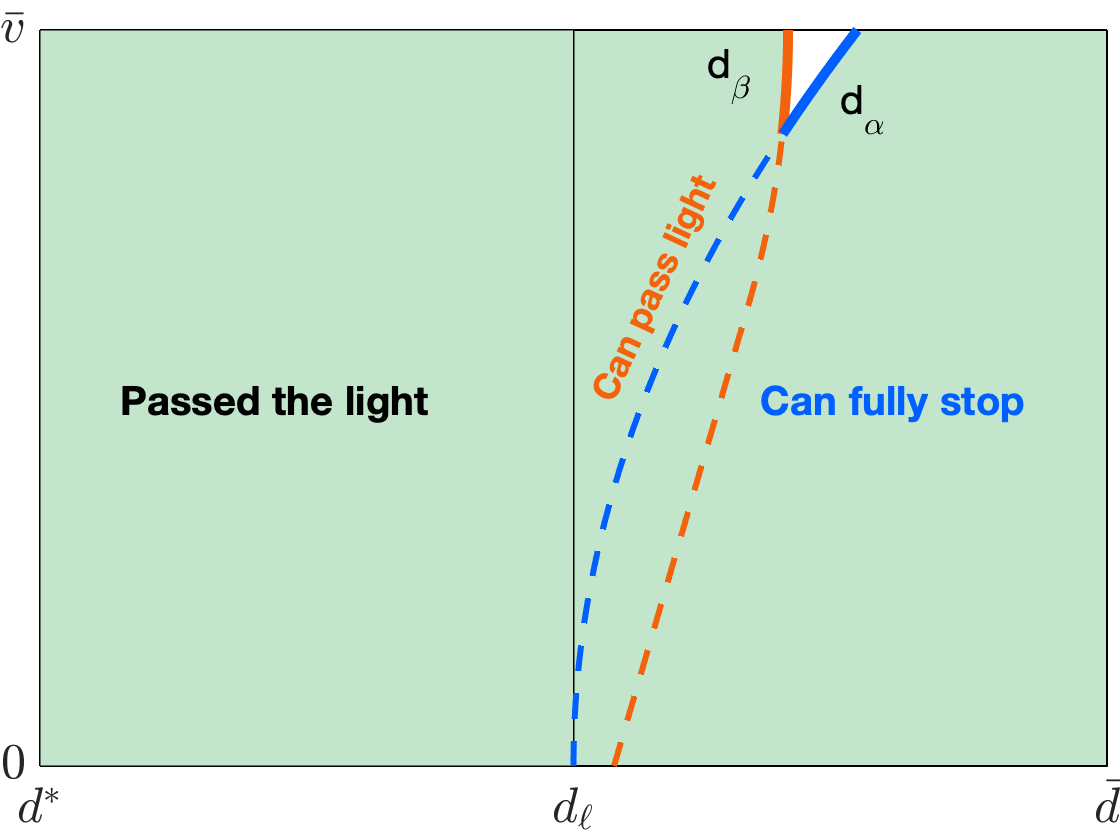}
\end{minipage}
\hspace{0.01cm}
\begin{minipage}[b]{0.45\linewidth}
\centering
\includegraphics[width=\textwidth]{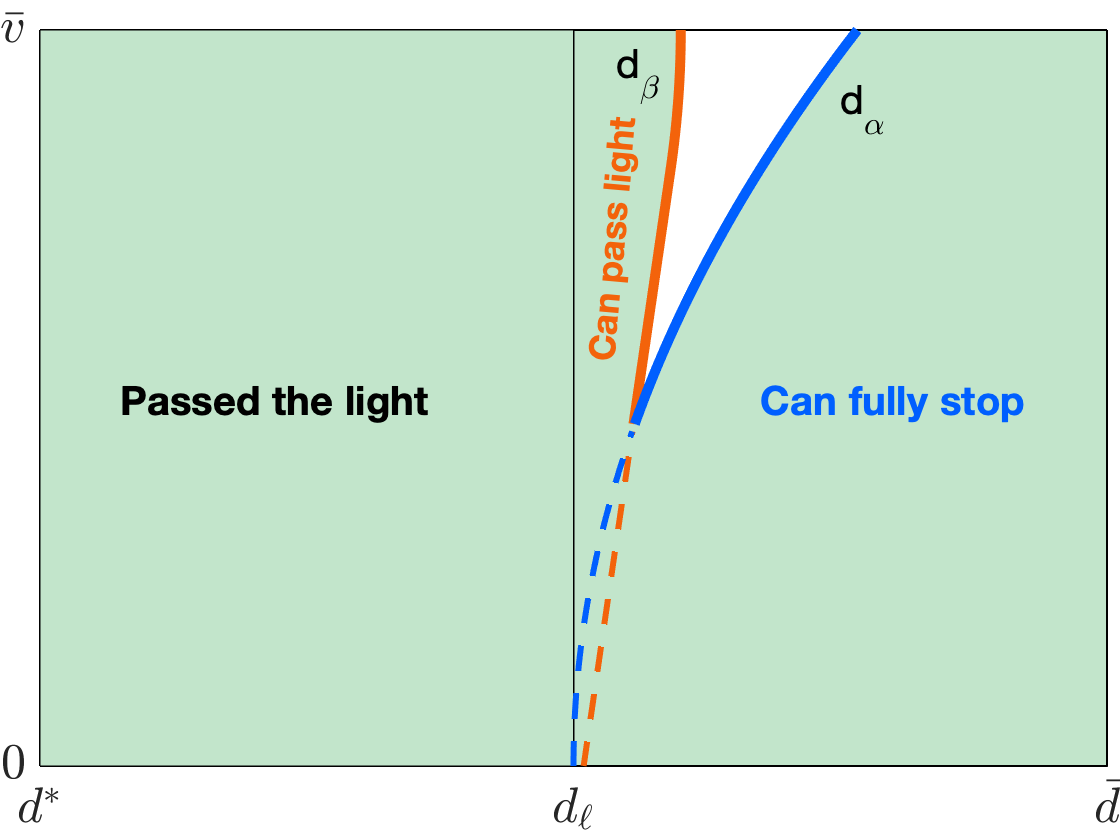}
\end{minipage}
\caption{
Two snapshots of the time-dependent piecewise boundary during phase $Y$. Allowed $(d,v)$ configurations are shown in green, and $\mathcal{I}\subY(t)$ is shown in white. 
For the starting configurations $(\dbeta(v,t), v)$, if speeding through the light before $t = \TR$ yields a lower cost than waiting for $T=\TG$, 
then the value function $u$ is discontinuous along the orange-dashed lines.
As $t \rightarrow \TR$, $\dbeta \rightarrow \dl$, and $\mathcal{I}\subY (t) \rightarrow \mathcal{I}\subR$.}\label{fig:yellow_phase_pw_bd_diagram}
\end{figure}

During the $Y$ phase, the car is allowed to have $d \leq d_\alpha(v)$ as long as it has enough time to enter the intersection by the time $\TR$.
The {\em largest} $d$ starting from which this is possible can be found by assuming that we use the maximum allowable acceleration (but without violating the speed limit) 
until reaching $\dl.$  
Denote the remaining yellow time as $\sr= \TR - t,$ the time that it would take to reach the speed limit as $\sbeta = \frac{\vbar - v}{\beta},$ 
and the lowest starting speed from which the speed limit would be reached before $\TR$ as $\vc = \vbar - \beta s.$  For fixed $v$ and $t \in [\TY, \, \TR),$ the largest allowed $d$ is then 
\begin{equation}\label{eqn:pw_max_accel_bd}
 \dbeta(v,t) =  \begin{cases}
 \dl + v\sr + \frac{\beta}{2}\srsq &  v \leq \vc;\\
 \dl + v \sbeta + \frac{\beta}{2}{\sbeta}^2 + \bar{v}(\sr - \sbeta), & v \geq \vc.
\end{cases} 
\end{equation}
The new disallowed region
$$\mathcal{I}\subY (t) \; = \; \left\{
(d,v) \mid \dbeta(v,t) \leq d \leq \dalpha(v), \, v \in [0, \bar{v}]
\right\}$$
is 
enforced for $t \in [\TY, \,  \TR)$ and
illustrated in Fig.\ \ref{fig:yellow_phase_pw_bd_diagram}. When $\dbeta(v,t) < \dalpha(v)$, the value function along $(\dbeta(v,t), v, t)$ is
\begin{equation}
\label{eq:c_beta}
C_{\beta} = (c_1 \beta + c_2 \beta^2)\shatbeta + c_3 \sr + u(\y(\TR), \TR)
\end{equation}
where $\shatbeta = \min\{\sbeta, \sr\}$. 
We assume that the yellow phase duration $\DY$ is long enough so that, starting from every $(d,v),$ it is possible to either enter the intersection before $\TR$ or 
come to a full stop before the intersection; i.e., $\mathcal{I}\subY (\TY) = \emptyset.$
Whenever $\dalpha(v) \leq \dbeta(v,t),$ the latter starting position allows either option, and if speeding toward the intersection yields a lower cost, then 
the value function will be discontinuous at $\left( \dbeta(v,t), \, v, t \right)$. 

\section{OPTIMAL DRIVING UNDER $\TY$ UNCERTAINTY}\label{section:uncertainty_models_paper}

Optimal planning in the initial green phase is quite similar if the duration ($D > 0$) and the starting time of that phase are fully known.
For a car that has been driving toward the intersection ever since this light turned green $\xi$ seconds before the current time $t=0$,
the natural planning horizon is 
$T = D - \xi.$ 
The light will then turn yellow, and the remaining cost-to-go will be 
$\delta(\x) = u(\x,\TY)$ already computed in Section \ref{section:problem_formulation_paper}.   
The new value function $w(\x,t)$ can be obtained by solving \eqref{eqn:hjb_timedep} with $t \in [0, \, T]$ 
and the terminal condition $w(\x,T) = \delta(\x).$
But in this section we focus on scenarios where $D$ is specified only probabilistically
while the durations of the subsequent phases ($\DY$ and $\DR)$ are known and fixed.
For a simple example, consider a city where all traffic light signals have green phase durations of either $D_1$ or $D_2 > D_1$ seconds,
and the fractions of signals in these respective categories are $p_1$ and $p_2=1-p_1$.  
Even though the driver is aware of the city-wide statistics, they may not know the category of this particular signal.
So,
the remaining green time 
is a random variable $\wTY$
with possible values $T_i = D_i - \xi$ and the corresponding probabilities $p_i, \; (i=1,2)$.
We are thus selecting a control $\ba(\cdot)$ to be used until the light turns yellow with the goal of minimizing
\begin{equation}
    \mathcal{C}(\x, t, a(\cdot)) = 
    \mathbb{E} \left\{
    \int^{\wTY}_t   \hspace*{-4mm} K \left( \y(r), a(r) \right) dr  +  \delta ( \y(\wTY) )
    \right\}.
\end{equation}
This is an example of ``initial uncertainty'' (\cite{QiDillonVlad_IU}),
with the problem becoming fully deterministic at the time $T_1$: the light will either turn yellow then
or 
we will immediately know that $\wTY = T_2.$
Thus, starting from any $t < T_1,$
\begin{align*}
    \mathcal{C}(\x, t, a(\cdot)) &=  
    \sum\limits_{i=1}^{2}
    p_i \left( 
    \int^{T_i}_t K \left( \y(r), a(r) \right) \, dr \; + \; \delta \left( \y(T_i) \right)
    \right)\\
    &= 
    \int^{T_1}_t K \left( \y(r), a(r) \right) \, dr \; + \; 
    p_1 \delta \left( \y(T_1) \right)
    \,+\,\\
    &
    (1-p_1) \left( 
    \int^{T_2}_{T_1} K \left( \y(r), a(r) \right) \, dr \; + \; \delta \left( \y(T_2) \right)
    \right). 
\end{align*}
This makes it natural to treat the planning on $[0,\, T_1)$ and $[T_1, \, T_2)$
as separate optimization problems, defining {\em a pair} of value functions $w^1(\x,t)$ and $w^2(\x,t)$
for the respective time intervals.
Each of these will solve the same PDE \eqref{eqn:hjb_timedep}, but with different terminal conditions
$$
w^1(x, T_1) \, = \, p_1 \delta(\x) + (1-p_1) w^2(\x,T_1);   \quad
w^2(\x,T_2) = \delta(\x).
$$

Generalizing this to $n$ possible green light durations $D_1 < \ldots < D_n$
with possible remaining green time intervals $T_i = D_i - \xi$ and  
$\PP\left(\wTY = T_i \right) = p_i$ for $i=1, \ldots, n,$
we introduce $n$ different value functions $w^i(\x,t),$ each defined on its own time interval
$[T_{i-1}, \, T_i),$ where for notational convenience we take $T_0 = 0.$
All of these satisfy \eqref{eqn:hjb_timedep} with the terminal conditions 
\begin{equation}
\label{eq:term_uncertain}
w^i \left(\x, T_{i} \right) \; =\;
\hat{p}_i \delta(\x) + (1-\hat{p}_i) w^{i+1} \left(\x,T_i \right),
\end{equation}
where $\hat{p}_i$ is the conditional probability 
$$
\hat{p}_i \; = \; \PP \left( \wTY = T_i \mid \wTY > T_{i-1} \right)
\; = \; 
p_i \, / \, \left( \sum\limits_{j=i}^n p_j \right).
$$ 
Since $\hat{p}_n = 1,$ the problem is deterministic on $[T_{n-1}, \, T_n).$
For $i < n,$ each $w^i$ depends only on the next $w^{i+1}$, and all value functions
can be found in a single sweep backward in time from $t=T_n$ to $t=T_0$.

\section{NUMERICAL IMPLEMENTATION}\label{section:numerical_implementation_paper}
\subsection{Numerical Methods for HJB Equation}
Overall, our numerical approach to the ``uncertain green phase duration'' problem
consists of the following three stages:
\begin{enumerate}
    \item Solve the stationary HJB \eqref{eqn:hjb_stationary} for the last phase $G$ .
    \item Solve HJB \eqref{eqn:hjb_timedep} for the deterministic phases $R$ and $Y$.
    \item Solve a sequence of HJBs \eqref{eqn:hjb_timedep} for the uncertain
    phase $I$. 
\end{enumerate}
The numerical results from each stage are used in determining the terminal conditions for the following stage. 
Our algorithms are based on a semi-Lagrangian (SL) discretization (\cite{Falcone_book}) on a Cartesian grid over a 
$(d,v,t)$ domain with $(N_d+1), (N_v+1),$ and $(N_t+1)$ gridpionts along the respective dimensions.
For a pre-specifed $N_v$, we select the discretization parameters
\begin{equation}
\label{eqn:deltas}
\Delta v = \frac{\bar{v}}{N_v}, \quad
\Delta t = \frac{\Delta v}{\max(\alpha, \beta)}, \quad 
\Delta d = \bar{v}\Delta t.
\end{equation}
The physical position, velocity, and time at a node $(i,j,k)$ are then
$(\x_{ij}, t_k) = (d_i, v_j, t_k) = \left( d^* + i\Delta d, \,  j\Delta v, \, k\Delta t \right)$ for
$i = 0, \dots, N_d$, $j = 0, \dots, N_v$, and $k = 0, \dots, N_t.$ 
The usual first-order SL discretization is
based on assuming that a fixed control $a$ is used starting from $\x_{ij}$
for a small time $\tau.$  The resulting 
new state $\tilde{\x}_{ij}^a = (\tilde{d}_i^a, \tilde{v}_j^a)$ is usually
approximated numerically (e.g., as $\tilde{\x}_{ij}^a \approx \x_{ij} + \tau \Bf(\x_{ij}, a)$),
but our simplified dynamics \eqref{eqn:frictionless_dynamics} allows computing $\tilde{\x}_{ij}^a$ analytically.
Moreover, \eqref{eqn:deltas} guarantees that $\tilde{d}_i^a \in [d_{i-1}, d_i]$ 
and $|\tilde{v}_j^a - v_j| \leq \Delta v$  for any $\tau \leq \Delta t,$ 
which makes it easier to interpolate the value function at $\tilde{\x}_{ij}^a.$
We also use the Golden Section Search (GSS) algorithm whenever we need to find the optimal $a$ numerically.

\subsubsection{Stage 1: Stationary HJB Solve}
In phase $G,$ there is never any incentive to decelerate and we 
only consider control values $a \in [0, \beta].$
The PDE \eqref{eqn:hjb_stationary} is stationary,
and the solution $\ug(\x_{ij})$ is approximated by a grid function
$Q_{ij}$ 
satisfying
\begin{align}
\nonumber
& Q^{}_{ij} \; = \; \min\limits_{a \in [0, \beta]} \left\{
\tau K(\x_{ij}, a)  + Q(\tilde{\x}_{ij}^a) \right\}, \quad i > 0, \, j < N_v;\\
\nonumber
& Q^{}_{i \mbox{\tiny N}^{}_v} = \; 
\tau K(\x_{i \mbox{\tiny N}^{}_v}, 0)  + Q(\tilde{\x}_{i \mbox{\tiny N}^{}_v}^0), \qquad \qquad \forall i;\\
\label{eqn:standard_stationary_SL}
& Q^{}_{0 j} \; = 0, \qquad \qquad \qquad \qquad \qquad \qquad \quad \; \forall j;
\end{align}
where the last two lines encode the need to stay under the speed limit and the boundary conditions respectively.
Since $\tilde{\x}_{ij}^a$
is usually not a gridpoint,
its value $Q $ is recovered via bilinear interpolation. I.e.,
$$
Q(\tilde{\x}_{ij}^a) \, = \, 
\gamma\tinySubs{1}(a) Q\tinySubs{i, j} + 
\gamma\tinySubs{2}(a) Q\tinySubs{i, j+1} + 
\gamma\tinySubs{3}(a) Q\tinySubs{i-1, j+1} + 
\gamma\tinySubs{4}(a) Q\tinySubs{i-1, j},
$$
with the nonnegative bilinear coefficients
adding up to 1 and $\gamma_1^{}(a) < 1$ for all $a$ values.
The coupled system \eqref{eqn:standard_stationary_SL} can be solved through {\em value iterations}, but
the dependence of the right hand side on $Q_{ij}$ slows down the convergence considerably even when using  a Gauss-Seidel relaxation.
Instead, we replace the first equation in \eqref{eqn:standard_stationary_SL} with an equivalent
\begin{equation}
\begin{split}
    Q_{ij} & =
    \min_{a \in [0, \beta]} \Bigg\{
   \frac{1}{1 - \gamma_{1}(a)} \Big(
   \tau K(\x_{ij}, a) + \gamma_2(a) Q_{i,j+1} \\
    & \qquad \qquad + \gamma_3(a)Q_{i-1, j+1} + \gamma_4(a)Q_{i-1,j} 
\Big) 
\Bigg\}.
\end{split}
\end{equation}
This effectively decouples the system, and we can now solve it in 
one
sweep, looping through $i=1, \ldots, N_d; \;  j  = N_v, ..., 0.$
\subsubsection{Stage 2: Deterministic $Y$-$R$ Phase HJB Solve}
The SL discretization for the time-dependent PDE \eqref{eqn:hjb_timedep} is similar,
but with several subtleties, which we describe below only briefly due to the space constraints.
The grid-approximation of the value function must satisfy
\begin{equation}
\label{eqn:sl_tHJB} 
    U^{k}_{ij} = \min_{a \in A_{ij}} \left\{K(\x_{ij}, a)\tau \, + \, U\left(\tilde{\x}_{ij}^a, t+\tau \right)\right\}
\end{equation}
for all $i,j$ and all $k < N_t$, which can be solved backwards in time from the terminal conditions $U^{N_t}_{ij} = Q_{ij}$.
The control set is
based on
the speed constraints: $A_{i0} = [0, \beta], \, A_{i \mbox{\tiny N}^{}_v} = [-\alpha, 0],$ and $A_{ij} =[-\alpha, \beta]$ for all other $j$ values.
On most of the domain we use $\tau = \Delta t,$ and the value at $(\tilde{\x}_{ij}^a, \, t+\tau)$ is obtained via bilinear interpolation in $(d,v)$.
But we reduce $\tau$ wherever it is needed to respect 
additional state constraints due to the traffic light, in which case we employ {\em cut-cells} in $(d,v)$ and additional linear interpolation in $t.$

In the red phase (with $t \in [\TR, \TG)$), the car has a disallowed region $\mathcal{I}\subR$.
For $\x_{ij}$ just right of the line $d = \dalpha(v)$, our $\tau$ is decreased adaptively, to ensure that $\tilde{d}_{ij}^a \geq \dalpha \left(\tilde{v}_{ij}^a \right).$
All grid cells intersected by this line are treated as cut-cells, with values of $U$ on the parabolic boundary 
computed by formula \eqref{eq:c_alpha}.  

Extending the solution into the yellow phase (for $t \in [\TY, \TR)$) presents two additional 
complications:
the time-dependent disallowed region $\mathcal{I}\subY(t)$
(see Fig. \ref{fig:yellow_phase_pw_bd_diagram}) and a possible discontinuity of $u(\x,t)$ when $d = \dbeta(v, t).$  
We adopt the following two-pass procedure to handle these challenges separately.
First, we solve \eqref{eqn:sl_tHJB} for all $d \geq \dalpha(v)$ exactly as we did in phase $R$; i.e., treating $d=\dalpha(v)$ as a boundary with the boundary conditions specified by  \eqref{eq:c_alpha}.  
This computes the best cost attainable by waiting out the red light.  But for some starting configurations, it is also possible to accelerate enough and cross the intersection before the light turns red. 
So, we then re-solve 
\eqref{eqn:sl_tHJB}
for all  $d \leq  \dbeta(v, t)$ without treating $d=\dalpha(v)$ as a boundary, with grid values updated only if they are smaller than those obtained in the first pass.  
To avoid interpolating across the discontinuity, in this second pass we use cut-cells just left of $d = \dbeta(v, t)$ 
with the boundary conditions specified by \eqref{eq:c_beta}.

\subsubsection{Stage 3: Uncertain Initial Phase $I$}
Given a set of possible times when the light could turn yellow $T_1, \ldots, T_n$,
each conditional value function $w^i(\x,t)$ satisfies PDE \eqref{eqn:hjb_timedep}
on $t \in [T_{i-1},T_i).$ 
Thus, we solve them sequentially (from $w^n$ to $w^1$) 
using the same SL-discretization \eqref{eqn:sl_tHJB} but  
with the terminal conditions \eqref{eq:term_uncertain}.
We note that the latter uses $\delta(\x) = u \left(\x, \TY \right)$ obtained at the end of Stage 2.
Once this function is computed, Stage 3 can be repeated as needed (for different probability distributions of $\wTY$) without 
repeating the computations in Stages 1 and 2.
We also note that the timestep in stage $I$ is always $\tau = \Delta t$ and cut-cells are not needed 
since 
$\mathcal{I}\subY(\TY) = \emptyset.$

\subsection{Optimal Trajectory Tracing}
Once the value function is computed, it is easy to find an approximately optimal trajectory
for every starting configuration $(d,v,t).$ 
For $t \in [\TY, \TG)$,
the optimal control value $a_{\star}$ is stored at every gridpoint while solving \eqref{eqn:sl_tHJB};
at all other points, we approximate $a_{\star}(\x,t)$ via tri-linear interpolation.
When tracing a trajectory, we use this approximately optimal feedback control to solve \eqref{eqn:frictionless_dynamics} 
until the time $\TG.$  At that point, we switch to the approximately optimal feedback control obtained from \eqref{eqn:standard_stationary_SL}
and continue until reaching the target $d^*.$

If starting in the initial green phase of uncertain duration from some $t \in [T_{i-1}, T_{i}),$
we again follow the feedback control obtained when solving for $w^i$.  
If the light turns yellow at the time $T_i$, we switch to a deterministic setting with the feedback control based on $u$. 
Otherwise, (i.e., if it stays green at $T_i$), we continue with $w^{i+1}.$

\subsection{Approximating the $\mathcal{J}_3$-Constrained Pareto Front}
Our running cost $K$ defined in \eqref{eqn:running_cost_ob_paper} can be viewed as a linear combination of constituent running costs
describing separately the fuel usage, discomfort from rapid speed changes, and time to target; i.e.,
$K_{1} = \gamma [a]_+$, $K_{2} = a^2$, and $K_{3} = 1$, with the coefficients $(c_1,c_2,c_3)$
reflecting the relative importance of these to the driver.  
Cumulative constituent costs $(\J_1, \J_2, \J_3)$ can be computed by integrating the 
corresponding $K_i$'s along the trajectory. 
Rational tradeoffs among these objectives can be explored by varying ratios of $c_i$'s.
To evaluate the fuel-discomfort tradeoffs among  trajectories
reaching the target in at most $T^{\dagger}$ seconds, we fix a ratio $\frac{c_1}{c_2}$ and find $c_3$ such that $\J_3 = T^{\dagger}$
for the $K$-optimal trajectory recovered from the HJB PDE.
Repeating this process for a range of $\frac{c_1}{c_2}$ values, we obtain a $J_3$-constrained {\em Pareto Front} for $(\J_1, \J_2).$
 
 \section{NUMERICAL EXPERIMENTS}\label{section:numerical_exp_paper}
 
For all experiments\footnote{
In the interest of computational reproducibility, we provide the full source code, additional figures, and movies for all examples at\\
\url{https://eikonal-equation.github.io/Traffic_Light_Uncertainty/}
}, we use realistic parameter values $\bar{v} = 45\, mph \approx 20.12\, m/s$, $-\alpha = -3.8\, m/s^2$, $\beta = 3.8\, m/s^2$, $\DY = 3 \, s$, $\DR = 60 \, s$.
We conduct all planning on a 200-meter-long road segment with the traffic light in the middle; i.e.,  $\bar{d} = 100 \, m$, $d_{\ell} = 0$, and $d^* = -100 \, m$. 
We use a $345 \times 181$ grid in $(d,v)$ space, with $\Delta t$ specified by \eqref{eqn:deltas}.

Mirroring the previous sections, we start with the planning in the later (deterministic) phases before considering the effects of uncertainty
in the initial phase $I$.

\begin{figure}[ht]
$
\begin{array}{cc}
\includegraphics[width=0.45\linewidth]{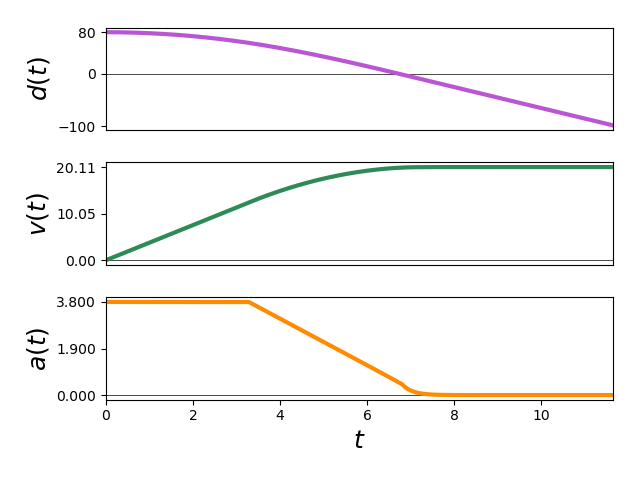} &
\includegraphics[width=0.45\linewidth]{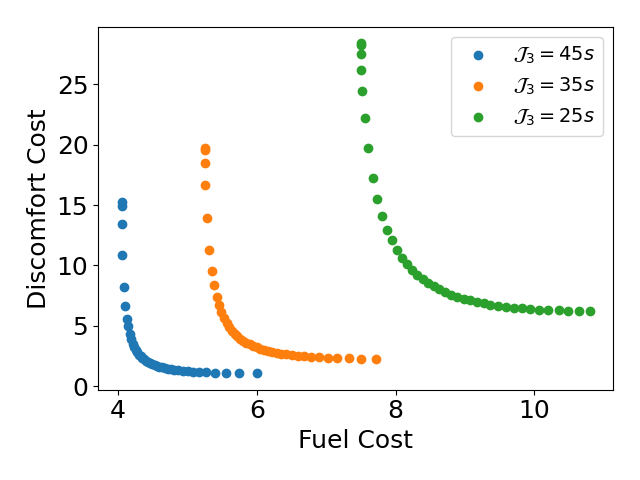}\\[-10pt]
\mbox{\footnotesize (A)} & \mbox{\footnotesize(B)}
\end{array}
$
\caption{Example 1 optimal trajectory and fuel-discomfort tradeoffs starting from $(d,v) = (80, 0)$. 
(A): Optimal $d(t)$, $v(t)$, and $a(t)$ when $(c_1,c_2,c_3) = (0.025, 0.025, 0.95)$.  The individual costs along the trajectory are $\mathcal{J}_1 = 20.11$, $\mathcal{J}_2 = 66.82$ 
and $\mathcal{J}_3 = 11.73$s. 
(B): $\mathcal{J}_3$-constrained Pareto fronts for $(\mathcal{J}_1, \mathcal{J}_2)$ with $\mathcal{J}_3$ = $25 s$, $35 s$, and $45s$. }\label{fig:green_only_pareto_and_traj}
\end{figure}

{\em  Example 1: Green phase $G$ only.}
The planning in this final green phase 
is particularly easy to interpret since there are no constraints other than the speed limit.
With our simplified dynamics, the method of characteristics used on PDE \eqref{eqn:hjb_stationary} shows that any optimal trajectory will have at most three stages:
(1) using the maximum acceleration $a=\beta$, (2) linearly decreasing acceleration until $a=0$, and (3) coasting until we reach the target.
We illustrate this for a specific starting $(d,v) = (80, 0)$ in Fig.\ \ref{fig:green_only_pareto_and_traj}A.  
Starting from the same $(d,v),$ Fig.\ \ref{fig:green_only_pareto_and_traj}B shows the rational  (fuel usage / acceleration discomfort) tradeoffs for 3 different constraint levels on the time-to-target.    

\begin{figure}[ht]
$
\begin{array}{cc}
\includegraphics[width=0.45\linewidth]{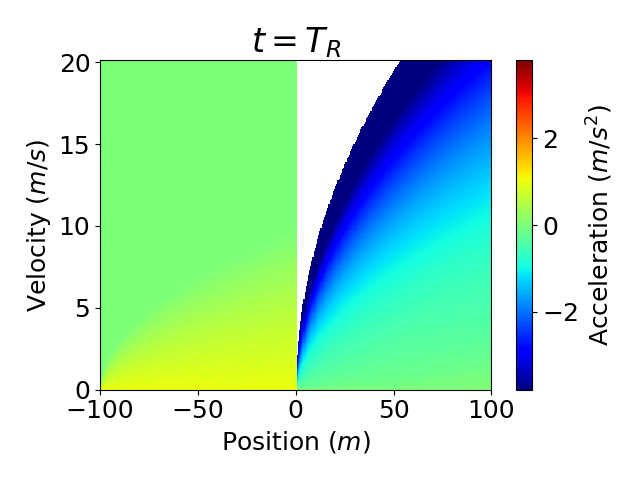} &
\includegraphics[width=0.45\linewidth]{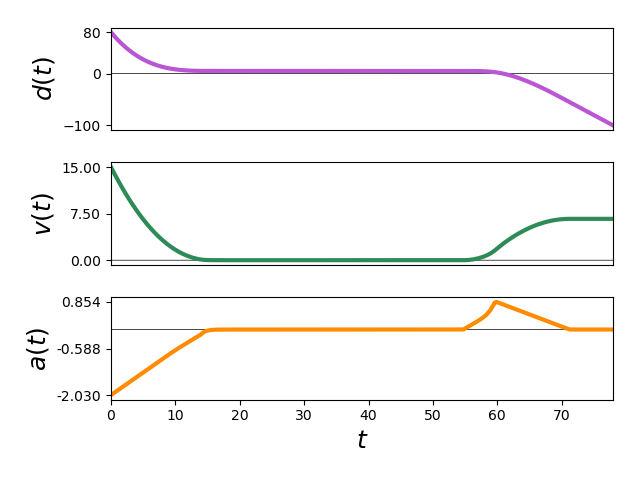}\\[-10pt]
\mbox{\footnotesize (A)} & \mbox{\footnotesize(B)}
\end{array}
$
\caption{Example 2 feedback controls and optimal trajectory when $(c_1, c_2, c_3) =( \frac{1}{3}, \frac{1}{3}, \frac{1}{3})$. (A): Feedback controls at $t = T\subR$ with $\mathcal{I}\subR$ shown in white. (B): Optimal $d(t), v(t)$ and $a(t)$ for starting $(d,v,t) = (80,15, T\subR)$. It is optimal for the driver to brake, stop, and wait until they can accelerate to arrive at $d_{\alpha}(v)$ at $T\subG$ with $v > 0$.}\label{fig:rg_trajectories}
\end{figure}

{\em Example 2: Phases $R$ and $G$. \hspace{1mm}}
Figure  \ref{fig:rg_trajectories} shows what happens when the driver starts planning as the light turns red at $t = \TR.$
Beyond the traffic light, for $d < \dl$, the optimal policy is exactly the same as in phase  $G.$
But on the other side, the car cannot cross the parabolic boundary $\dalpha(v) = d$ until $\TG = \TR + 60s$.
This requires aggressive braking to the right of this parabola, though the need for harsh braking diminishes as $t \rightarrow T\subG$ and $d \rightarrow \bar{d}$ since it is optimal to reach $\dalpha(v) = d$ with non-zero speed just as the light turns green.

\begin{figure}[ht]
$
\begin{array}{cc}
\includegraphics[width=0.45\linewidth]{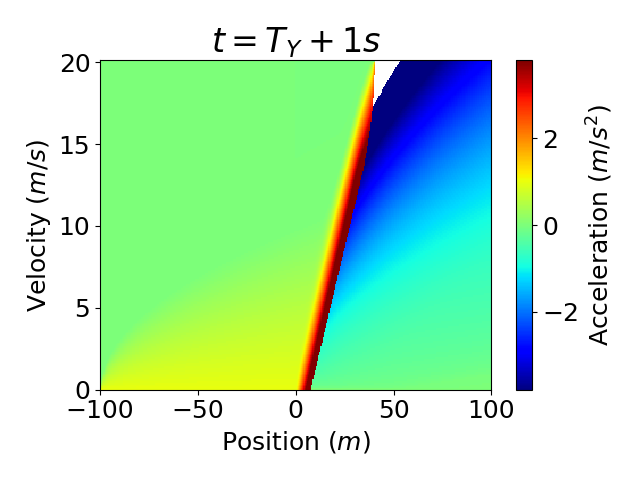} &
\includegraphics[width=0.45\linewidth]{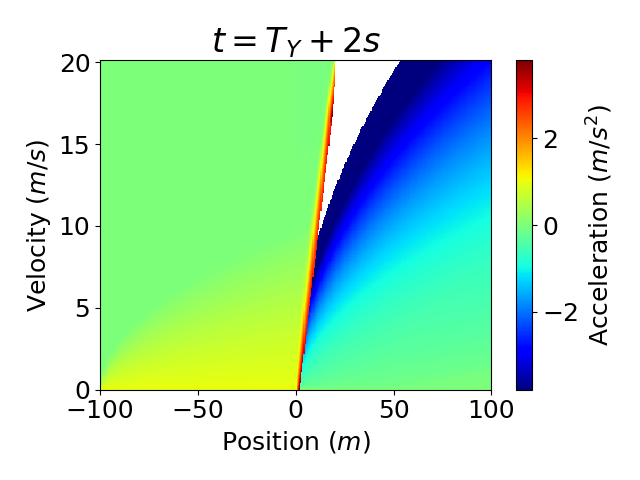}\\[-10pt]
\mbox{\footnotesize (A)} & \mbox{\footnotesize(B)}
\end{array}
$
\caption{Example 3 feedback controls for $(c_1, c_2, c_3) =( \frac{1}{3}, \frac{1}{3}, \frac{1}{3})$ at (A): $t = \TY +1 \, s$ and (B): $t = \TY + 2\, s$. 
Both plots are the control-heatmap equivalents to the piecewise boundary schematics shown in Fig. \ref{fig:yellow_phase_pw_bd_diagram}.}\label{fig:yrg_oc}
\end{figure}

\begin{figure}[ht]
$
\begin{array}{cc}
\includegraphics[width=0.45\linewidth]{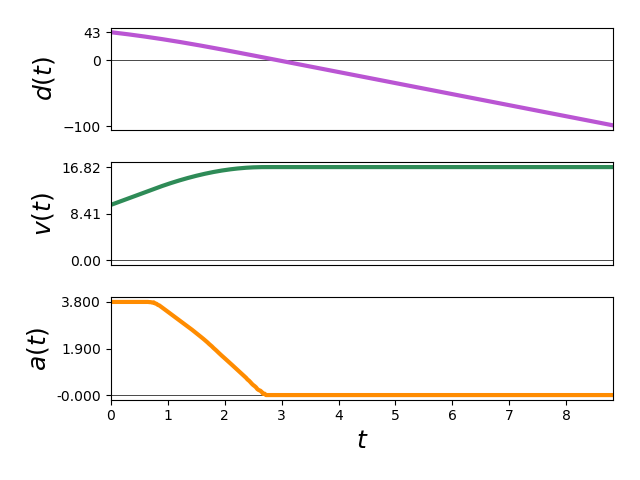} &
\includegraphics[width=0.45\linewidth]{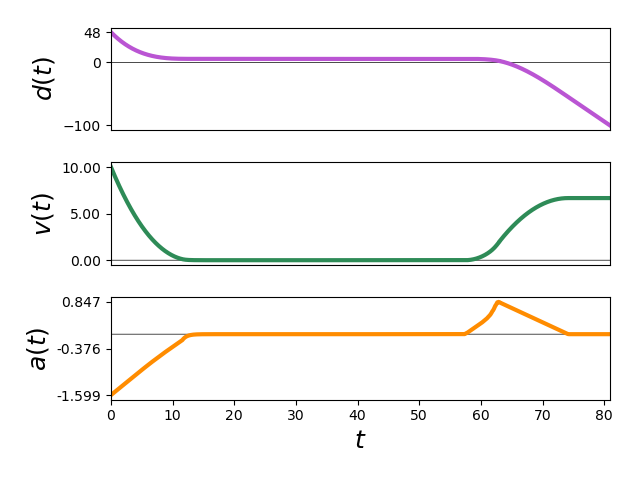}\\[-10pt]
\mbox{\footnotesize (A)} & \mbox{\footnotesize(B)}
\end{array}
$
\caption{Example 3 optimal trajectories starting from (A): $(d,v,t) =  (43, 10, T\subY)$ and (B): $(d,v,t) =  (48, 10, T\subY).$
The former yields $(\mathcal{J}_1, \mathcal{J}_2, \mathcal{J}_3) = (6.82, 21.32, 8.9)$ and the overall cost $\mathcal{J} = 12.35.$
The latter results in $(\mathcal{J}_1, \mathcal{J}_2, \mathcal{J}_3) = (6.65, 14.16, 80.95)$ and a much larger $\mathcal{J} = 33.94.$
}\label{fig:yrg_traj}
\end{figure}

{\em Example 3: Deterministic phases $Y, R,$ and $G.$}
In the yellow phase, some of the starting configurations (with $d \leq \dbeta(v,t)$) allow beating the red light if the driver is ready to accelerate aggressively.
This leads to a time-dependent disallowed set $\mathcal{I}_Y(t)$ (see Fig. \ref{fig:yrg_oc}) and a possible discontinuity in the value function when $d = \dbeta(v,t).$
Fig. \ref{fig:yrg_traj} shows that even a slight shift in a starting position over this discontinuity will result in a drastically different optimal trajectory and cumulative cost.   

\begin{figure}[ht]
$
\begin{array}{cc}
\includegraphics[width=0.45\linewidth]{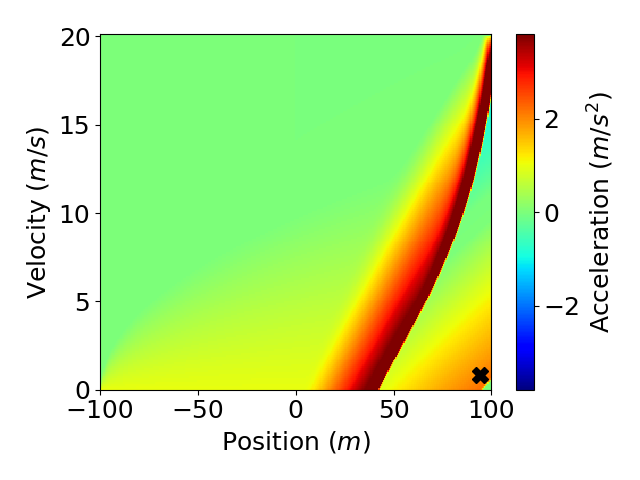} &
\includegraphics[width=0.45\linewidth]{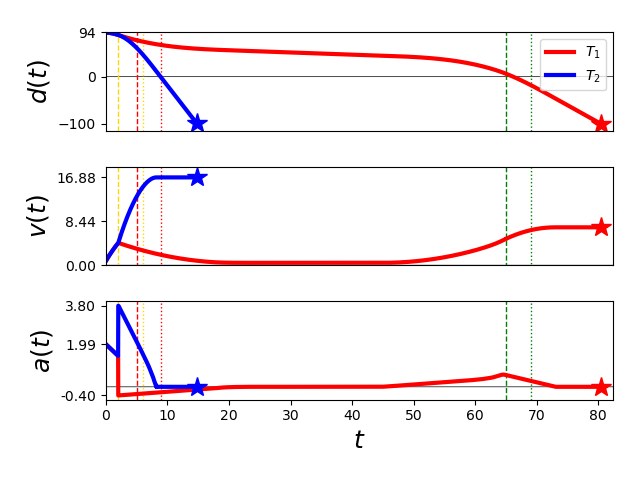}\\[-10pt]
\mbox{\footnotesize (A)} & \mbox{\footnotesize(B)}
\end{array}
$
\caption{Example 4 feedback controls and trajectories for starting $(d,v) = (94,0.85)$ when $(p_1, p_2) = (0.5,0.5)$ and $(c_1, c_2, c_3) =( \frac{1}{3}, \frac{1}{3}, \frac{1}{3})$. (A): Feedback controls at $t = 0$. The ``X" indicates the vehicle's starting point. (B): Optimal $d(t)$, $v(t)$ and $a(t)$. Vertical lines corresponding to the possible turning yellow, red, and green times are shown in their respective colors.}\label{fig:model1_2l}
\end{figure}

{\em Example 4: Uncertain yellow change, $\wTY \in \{T_1, T_2\}.$}\\
We now start planning at the time $t=0$ in the initial green phase $I,$ with the remaining time until yellow 
$\wTY$ taking values $T_1 = 2 \, s$ or $T_2 = 6 \, s$ with equal probability ($p_1= p_2 = 1/2$).
We first assume that the driver values all three objectives equally ($c_1=c_2=c_3 = \frac{1}{3}$);
the resulting optimal feedback control for $t=0$ is shown in Fig. \ref{fig:model1_2l}A.
We focus on the starting configuration $(d,v) = (94,0.85),$ from which it is possible to beat the traffic light
if $\wTY = T_2$ but not if  $\wTY = T_1.$  The optimal trajectory (shown in Fig. \ref{fig:model1_2l}B) 
branches at the time $t=T_1$, when we find out the true value of $\wTY$.
We note that the ``optimal under uncertainty'' control used for $t \in [0,T_1)$ would not be optimal for either deterministic scenario:
we would accelerate far less (if at all) with $\TY = T_1$ and far more with $\TY=T_2.$

\begin{figure}[ht]
$
\begin{array}{cc}
\includegraphics[width=0.45\linewidth]{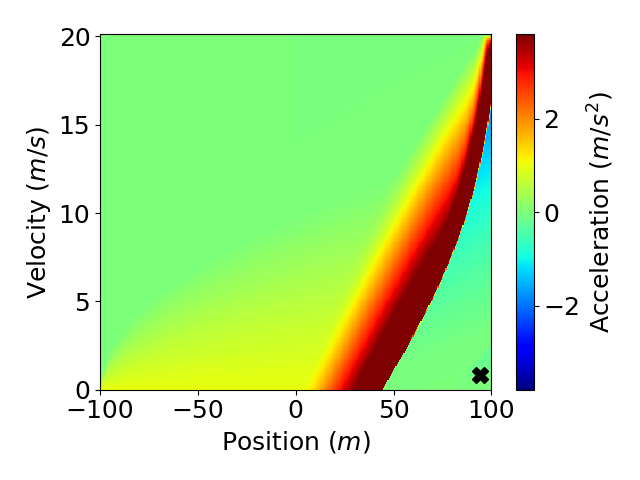} &
\includegraphics[width=0.45\linewidth]{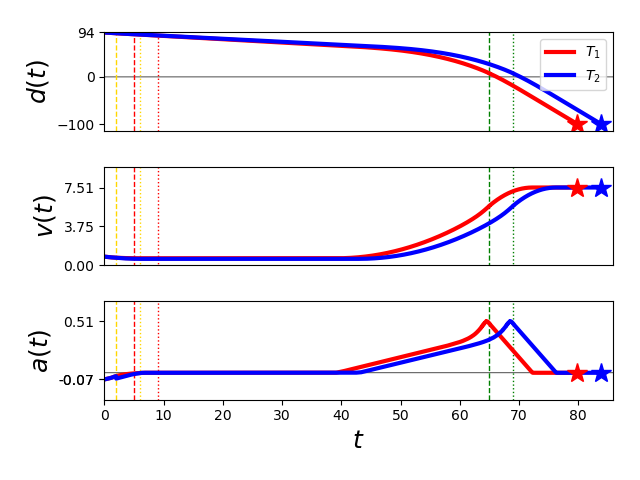}\\[-10pt]
\mbox{\footnotesize (A)} & \mbox{\footnotesize(B)}
\end{array}
$
\caption{Example 4 feedback controls and trajectories for starting $(d,v) = (94,0.85)$ when $(p_1, p_2) = (0.95, 0.05)$, and $(c_1, c_2, c_3) =( \frac{1}{3}, \frac{1}{3}, \frac{1}{3})$. (A): Feedback controls at $t = 0$.  The ``X" indicates the vehicle's starting point. (B): Optimal $d(t)$, $v(t)$ and $a(t)$. Vertical lines corresponding to the possible turning yellow, red, and green times are shown in their respective colors.}\label{fig:model1_2l_t1likely}
\end{figure}

\begin{figure}[ht]
$
\begin{array}{cc}
\includegraphics[width=0.45\linewidth]{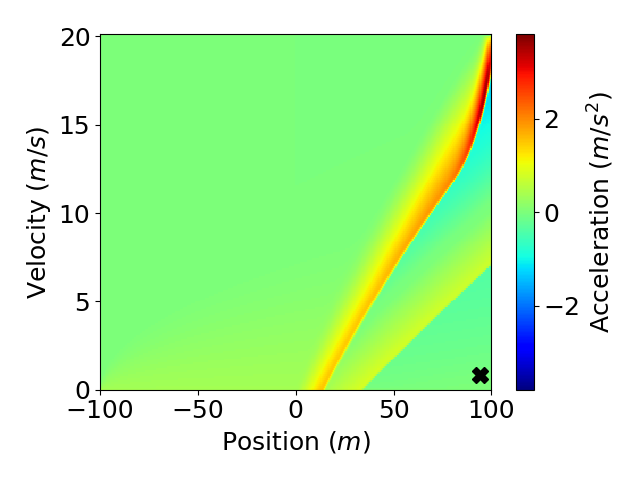} &
\includegraphics[width=0.45\linewidth]{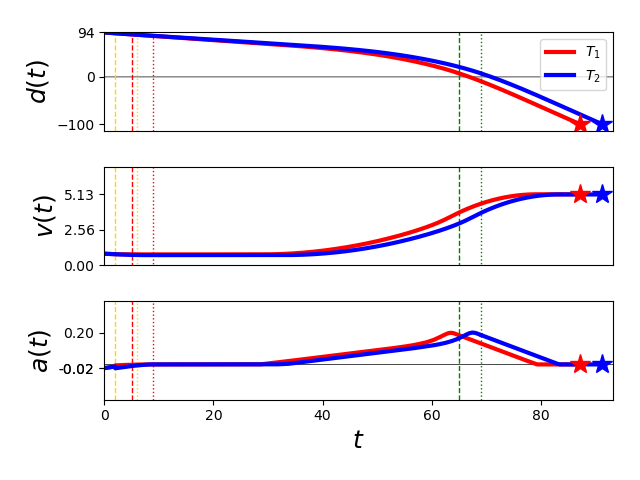}\\[-10pt]
\mbox{\footnotesize (A)} & \mbox{\footnotesize(B)}
\end{array}
$
\caption{Example 4 feedback controls and trajectories for starting $(d,v) = (94,0.85)$ when $(p_1, p_2) = (0.5, 0.5)$, and $(c_1, c_2, c_3) = (0.15, 0.75, 0.1)$. (A): Feedback controls at $t = 0$. The ``X" indicates the vehicle's starting point. (B): Optimal $d(t)$, $v(t)$ and $a(t)$. Vertical lines corresponding to the possible turning yellow, red, and green times are shown in their respective colors.}\label{fig:model1_2l_c2greatest}
\end{figure}

But just because one can beat the traffic light with $\wTY = T_2,$ it does not mean that it is always optimal to do so.
Indeed, if this light duration is sufficiently unlikely, the initial acceleration on $[0,T_1)$ to preserve both options is no longer worthwhile.
We demonstrate this for $(p_1, p_2) = (0.95, 0.05)$ in Fig.  \ref{fig:model1_2l_t1likely}.  (Note that there is still branching at $t=T_1$ since at that point we discover the time
$(\wTY+ \DY + \DR)$ by which we need to reach the parabolic boundary as the light turns green again.)  Not surprisingly, a similar decision not to rush can also result from a difference 
in driver's priorities; e.g., a high enough $c_2$ will make rapid acceleration unattractive even if $T_2$ is fairly likely.  In Fig. \ref{fig:model1_2l_c2greatest}
we demonstrate this for $(p_1, p_2) = (0.5, 0.5)$ and $(c_1, c_2, c_3) =( 0.15, 0.75, 0.1)$.

\begin{figure}[t]
$
\begin{array}{cc}
\includegraphics[width=0.45\linewidth]{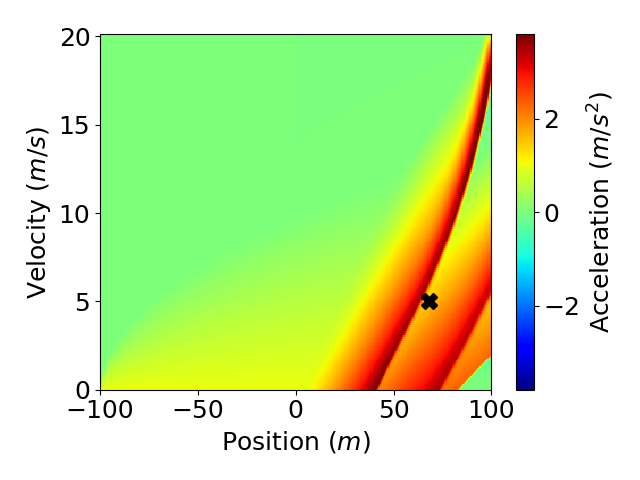} &
\includegraphics[width=0.45\linewidth]{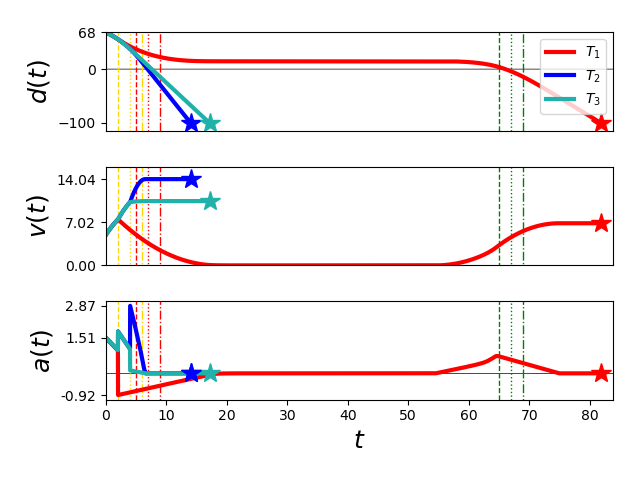}\\[-10pt]
\mbox{\footnotesize (A)} & \mbox{\footnotesize(B)}
\end{array}
$
\caption{Example 5 feedback controls and trajectories for starting $(d,v) = (68,5)$ when $(p_1, p_2, p_3) = (0.25, 0.25, 0.5)$, and $(c_1, c_2, c_3) =( \frac{1}{3}, \frac{1}{3}, \frac{1}{3})$. (A): Feedback controls at $t = 0$. The ``X" indicates the vehicle's starting point. (B): Optimal $d(t)$, $v(t)$, and $a(t)$. The $(T_1, T_2, T_3)$ and their corresponding $\TR$ and $\TG$ are marked by the yellow, red, and green vertical lines respectively.}\label{fig:model1_3l}
\end{figure}

{\em Example 5: Uncertain yellow change, $\wTY \in \{T_1, T_2, T_3\}.$}\\ 
 In our final example, we use $(T_1, T_2, T_3) = (2 \, s, 4 \, s, 6 \,s)$ with $(p_1, p_2, p_3) = (0.25, 0.25, 0.5)$ and $c_1 = c_2 = c_3= \frac{1}{3}$ 
 starting from $(d,v) = (68,5).$  The optimal control (shown in Fig. \ref{fig:model1_3l}) has three branches:  if the light stays green at $t=T_1,$
 it becomes clear that we can beat the red light, but how much acceleration will be needed to do so will be revealed 
 at $t=T_2$.

 \section{CONCLUSIONS}

We provide a framework for determining optimal driving strategies in the face of initial traffic light uncertainty and competing optimization objectives.  
While our basic setup is intentionally simple, the proposed optimization under uncertainty approach is much broader.
We hope that it will be useful in modifying prior models with detailed vehicle dynamics (\cite{sun2020optimal}), route selection on complex road networks with many signalized intersections (\cite{mahler2014optimal}), and game-theoretic traffic effects due to independent decision making of multiple drivers (\cite{huang2019game}).
A similar approach will also be useful in treating other uncertainty models; e.g., due to pedestrian-actuated signal timing changes (\cite{eteifa2021predicting}).

\textbf{Acknowledgements: } 
The authors are grateful to A.~Nellis, J.~van Hook, and N.~Do for their preliminary work on the Phase $R$ problem during the summer 2018 REU Program.
 
\bibliography{bibl}

\end{document}